\input amstex
\input amsppt.sty
\input label.def
\input degt.def
\input elliptic.def


{\catcode`\@11
\gdef\proclaimfont@{\sl}}

\newtheorem{Remark}\Remark\thm\endAmSdef

\def\dash{\item"\hfill--\hfill"}
\def\Dashes{\widestnumber\item{--}\roster}
\def\endDashes{\endroster}

\def\ie{\emph{i.e.}}
\def\eg{\emph{e.g.}}
\def\cf{\emph{cf}}
\def\via{\emph{via}}
\def\etc{\emph{etc}}

\def\I{\roman{I}}
\def\If{\I_{\operatorname{F}}}
\def\Ib{\I_{\operatorname{B}}}
\def\II{\roman{II}}
\def\inmath#1{\ifmmode#1\else$#1$\fi}
\def\1{\inmath1}
\def\2{\inmath2}
\def\3{\inmath3}
\def\n#1{\inmath{\bar#1}}

\def\BR{B_{\R}}
\def\CR{C_{\R}}

\def\XR{X_{\R}}

\def\SigmaR{\Sigma_{\R}}
\def\subIm{_{\Impart}}

\def\BIm{B\subIm}
\def\CIm{C\subIm}
\def\tE{\tilde E}
\def\tC{\tilde C}

\def\Sk{\operatorname{Sk}}
\def\Skdir{\Sk_{\fam0 dir}}
\def\Skud{\Sk_{\fam0 ud}}
\def\bbase{\bar\base}

\def\vector#1{\!\vec{\,#1}}


\topmatter

\title
Real trigonal curves and\\
real elliptic surfaces
of type~$\I$
\endtitle

\author
Alex Degtyarev, Ilia Itenberg, and Victor Zvonilov
\endauthor

\thanks
The second author is partially funded by
the ANR-09-BLAN-0039-01 grant
of {\it Agence Nationale de la Recherche}
and is a member of
FRG: Collaborative Research: Mirror Symmetry \& Tropical Geometry
(Award No. $0854989$).
\endthanks

\address
Bilkent University\endgraf\nobreak
Department of Mathematics\endgraf\nobreak
06800 Ankara, Turkey
\endaddress

\email
degt\@fen.bilkent.edu.tr
\endemail

\address
Universit\'e de Strasbourg,
IRMA and Institut Universitaire de France
\endgraf\nobreak
7 rue Ren\'e Descartes\endgraf\nobreak
67084 Strasbourg Cedex, France
\endaddress

\email
ilia.itenberg\@math.unistra.fr
\endemail

\address
Syktyvkar State University\endgraf\nobreak
Department of Mathematics\endgraf\nobreak
167001 Syktyvkar, Russia
\endaddress

\email
zvonilov\@syktsu.ru
\endemail


\abstract
We study real trigonal curves and elliptic surfaces of type~$\I$
(over a base of an arbitrary genus)
and their fiberwise equivariant deformations.
The principal tool is a real version of Grothendieck's \emph{dessins
d'enfants}.
We give a description of maximally inflected
trigonal curves of
type~$\I$ in terms of the combinatorics of sufficiently simple
graphs and,
in the case of the rational base,
obtain a complete classification of such curves.
As a consequence, these results lead to
conclusions concerning real Jacobian
elliptic surfaces of type~$\I$ with all singular fibers real.
\endabstract



\endtopmatter

\document

\section{Introduction}
This paper is a continuation of~\cite{DIK} and~\cite{Zvonilov},
where the authors, partially in collaboration with V.~Kharlamov,
have obtained a
complete deformation classification of the so called
real trigonal $M$- and
$(M-1)$-curves in geometrically ruled surfaces (see
Subsection~\ref{real-trigonal} for the precise settings).
Recall that a real algebraic or analytic variety is called an
$M$-{\it variety}
if it is maximal in the sense of the
Smith--Thom
inequality. A generalization of the notion of $M$-curves are the
curves of type~$\I$, \ie, those whose real part separates the set
of complex points. All $M$-curves are indeed of type~$\I$.
In the case of trigonal curves,
there is another
(almost) generalization:
one can consider a curve~$C$ such that
all critical points of the restriction to~$C$ of the ruling
of the ambient surface are real. We call such curves
\emph{maximally inflected}. According to~\cite{DIK}, for
a trigonal
$M$-curve~$C$, all but at most four critical points are real and,
moreover, the curve can be deformed to an essentially unique
maximally inflected one. In the present paper, we make an attempt
to adapt
the techniques used in~\cite{DIK} to maximally inflected
trigonal curves of type~$\I$, obtaining a complete classification
in the case of the rational base.
As usual, \cf\. \eg~\cite{DIK},~\cite{Bihan.Mangolte}, \etc.,
using the computation
of the real version of the Tate--Shafarevich group found
in~\cite{DIK}, one can extend these results, almost literally,
to real elliptic
surfaces.

Throughout the paper, all varieties are over~$\C$
(possibly, with a real structure)
and nonsingular.

\subsection{Principal results}
As in~\cite{DIK} and~\cite{Zvonilov},
the principal tool used in the paper
is the notion of \emph{dessin}, see
Section~\ref{rtc},
which is a real
version of Grothendieck's \emph{dessin d'enfants} of the
functional
$j$-invariant of the curve;
this concept was originally suggested by
S.~Orevkov~\cite{Orevkov} and then developed in~\cite{DIK}, where
the study of deformation classes of real trigonal curves was
reduced to that of dessins, see Proposition~\ref{equiv.curves}.
In the
case of maximally inflected curves of type~$\I$, we manage
to
simplify
the rather
complicated combinatorial structure of dessins to somewhat
smaller
graphs, which we call \emph{skeletons}, see
Section~\ref{S.skeletons}. One of our principal results is
Theorem~\ref{cor.Sk}, which
establishes a one-to-one correspondence between the equivariant
fiberwise deformation classes of maximally inflected trigonal
curves of type~$\I$ and certain equivalence classes of skeletons.

In the case of the rational base (\ie, for curves in
rational ruled surfaces),
skeletons can be regarded as unions of chords in the disk and their
equivalence takes an especially simple form. We use
Theorem~\ref{cor.Sk} and show that, in this case,
a trigonal curve of type~$\I$ is essentially determined by its
real part. More precisely, we prove the
following two statements (see
Subsections~\ref{proof.existence}
and~\ref{proof.main}, respectively).

\theorem\label{th.existence}
A maximally inflected
real trigonal curve~$C$ in a real rational
geometrically ruled surface~$\Sigma$ is of type~$\I$
if and only if its real part~$C_\R$ admits a quasi-complex
orientation, see Definition~\ref{def.qc.orientations}.
\endtheorem

\theorem\label{th.main} Let $\Sigma \to \Cp1$ be a real rational
geometrically ruled surface, and $C', C''\subset \Sigma$ two
maximally inflected real trigonal curves of type~$\I$. Then, any
fiberwise auto-homeomorphism of $\Sigma_\R$ isotopic to identity and
taking $C'_\R$ to $C''_\R$ is realized
by a fiberwise equivariant deformation
\rom(see~\ref{deformation}\rom)
of the
curves.
\endtheorem

An attempt
of a constructive description
of the real parts realized by maximally inflected type~$\I$
trigonal curves over the rational base is made in
Subsection~\ref{s.blocks}.

Note that, in the literature, there is a great deal of various
definitions of type~$\I$, especially in the case of surfaces:
usually, one requires that the real part of the variety should
realize $\bmod2$
a certain `universal' class in the homology of
the complexification. For trigonal curves and elliptic surfaces,
we introduce the notions of type~$\Ib$ and~$\If$, respectively,
see Subsections~\ref{s.Ib} and~\ref{s.If}.
While sharing most properties of the classical type~$\I$,
these notions are particularly well suited
for real trigonal curves and elliptic surfaces, extending the
concept of type~$\I$ to the case of non-separating base.

When working with trigonal curves and elliptic surfaces, the
ruling is regarded as part of the structure and, hence, the natural
equivalence relation is fiberwise equivariant deformation. It is
this equivalence that is dealt with in the bulk of the paper.
However, in general in topology of real algebraic varieties, it is
more common to consider a weaker relation, the so called
\emph{rigid isotopy}, which does not take the ruling into account.
A brief discussion of rigid isotopies of real trigonal curves is
found in Appendix~\ref{ap.A}. We prove
Theorem~\ref{max.inflected}, that states that any non-hyperbolic
(see~\ref{real-trigonal})
curve of
type~$\Ib$ is rigidly isotopic to a maximally inflected one (and,
in particular, the assumption that the curve should be maximally
inflected in the other statements is not very restrictive). Note
though, that this assertion is indeed specific for type~$\I$, see
Example~\ref{example}.

\subsection{Contents of the paper}
Sections~\ref{S.curves} and~\ref{rtc} are introductory:
we recall a few notions and facts
related to topology of real trigonal curves and their dessins,
respectively.
The concepts of type~$\Ib$ for trigonal curves and type~$\If$ for
elliptic surfaces are introduced and the relation between them is
discussed in Section~\ref{S.curves}.
In Section~\ref{curves.typeIb}, we study properties of dessins
specific to curves of type~$\Ib$, first in general, and then in
the maximally inflected case.
The heart of the paper is Section~\ref{S.skeletons}: we introduce
skeletons, define their equivalence, and prove
Theorem~\ref{cor.Sk}. Section~\ref{S.rational} deals with the case
of the rational base: we prove Theorems~\ref{th.existence}
and~\ref{th.main} and introduce \emph{blocks}, which are the
`elementary pieces' constituting the dessin of
any maximally inflected curve of
type~$\I$ over~$\Cp1$.
Finally, in Appendix~\ref{ap.A} we digress to not necessarily
fiberwise equivariant deformations of real trigonal curves and
show that, by such a deformation, all singular fibers of a
non-hyperbolic curve of type~$\I$ can be made real.

\subsection{Acknowledgments}
We are grateful to the \emph{Mathematisches
Forschungsinstitut Oberwolfach} and its RiP program for the hospitality
and excellent working conditions which helped us to complete this project.
A part of the work was done during the first author's visits
to \emph{Universit\'e de Strasbourg}.

\section{Trigonal curves and elliptic surfaces}\label{S.curves}

In this section, we recall a few basic notions and facts
related to topology of real trigonal curves,
introduce curves of type~$\Ib$ and elliptic surfaces of~$\If$,
and discuss the relation between these objects.

\subsection{Real trigonal curves}\label{real-trigonal}
Let $\pi\:\Sigma\to B$ be a geometrically ruled surface over a
base~$B$
and with the exceptional section~$E$, $E^2=-d<0$.
The fibers of the ruling~$\pi$ are often called
\emph{vertical}, \eg, we speak about vertical tangents, vertical
flexes \etc. We identify~$E$ and~$B$ \via\ the restriction of~$\pi$.
Denote by $e,f\in H_2(\Sigma)$ the classes realized by~$E$ and a
generic fiber~$F$, respectively.

A \emph{trigonal curve} on~$\Sigma$ is
a reduced curve $C\subset\Sigma$ disjoint from~$E$ and such that the
restriction $\pi_C\:C\to B$ of~$\pi$
has
degree three.
One has
$[C]=3e+3df\in H_2(\Sigma)$. Given~$C$,
we denote by $B^\circ$ the
complement
in~$B$
of
the critical locus of~$\pi_C$.

Given
a trigonal curve
$C\subset\Sigma$, the fiberwise center of gravity of the three
points of~$C$ (viewed as points in the affine fiber of
$\Sigma\sminus E$) defines an additional section~$0$ of~$\Sigma$;
thus, a necessary condition for $\Sigma$ to contain a trigonal curve
is that the $2$-bundle
whose projectivization is~$\Sigma$ splits.

Recall that a \emph{real
variety}
is a complex algebraic (analytic) variety~$V$ equipped with
an anti-holomorphic involution $c = c_V\: V \to V$;
the latter involution is called a \emph{real structure}
on~$V$.
The fixed point set $V_\R = \Fix c$ is called the \emph{real part}
of~$V$. A regular morphism $f\:V \to W$ of two real varieties
is called \emph{real} or \emph{equivariant}
if it commutes with the real structures, \ie, one has
$f\circ c_V=c_W\circ f$.

Let $\pi\: \Sigma \to B$ as above be real.
Throughout the paper we assume that $B_\R\ne\varnothing$.
The exceptional section $E \subset \Sigma$ is also real
and~$\pi$ establishes a bijection between
the connected components~$\Sigma_i$ of $\Sigma_\R$
and the connected components~$B_i$ of $B_\R$.
All restrictions $\pi_i\: \Sigma_i \to B_i$ are $S^1$-bundles,
not necessarily orientable. The sum $\sum w_1(\pi_i)[B_i]$
equals $d \bmod 2$.

Let $C \subset \Sigma$ be a nonsingular real trigonal curve.
The connected components of $C_\R$ split into
groups $C_i = C_\R \cap \pi^{-1}_C(B_i)$.
Given~$C$, a component~$B_i$ (and the group~$C_i$)
is called \emph{hyperbolic} (\emph{anti-hyperbolic})
if the restriction $C_i \to B_i$ of~$\pi$ is three-to-one
(respectively, one-to-one). The trigonal curve~$C$
is called \emph{hyperbolic} if all its groups are hyperbolic.

Each non-hyperbolic group~$C_i$ has a unique
\emph{long component} $l_i$, characterized by the fact
that the restriction $l_i \to B_i$ of~$\pi$
is of degree~$\pm 1$. All other components of~$C_i$
are called \emph{ovals}; they are mapped to~$B_i$
with degree~$0$. Let $Z_i \subset B_i$ be the set of points
with more than one preimage in $C_i$. Then, each oval
projects to a whole component of~$Z_i$, which is also called
an oval. The other components of $Z_i$,
as well as
their preimages
in~$l_i$, are called \emph{zigzags}.

A trigonal curve $C \subset \Sigma$
is called \emph{almost generic}
if it is nonsingular and all critical points
of the restriction~$\pi_C$
are simple; in other words, $C$ is almost generic if all its
singular fibers are of Kodaira type~$\I_1$
(or $\tilde\bA_0^*$ in the alternative notation).
A real trigonal curve~$C$ is called \emph{maximally inflected}
if it is almost generic and
all critical points of the restriction~$\pi_C$
are real.

\subsection{Deformations}\label{deformation}
Throughout this paper, by a \emph{deformation} of a trigonal curve
$C\subset\Sigma$ over~$B$ we mean a deformation of
the pair $(\pi\: \Sigma \to B, C)$
in the sense of Kodaira--Spencer.
It is worth emphasizing that the complex structure
on~$B$ and~$\Sigma$ is not assumed fixed;
it is also subject to deformation. (In the correspondence between
trigonal curves and dessins, see Proposition~\ref{equiv.curves}
below, the complex structure on
the base~$B$ is recovered using the Riemann
existence theorem.)
A deformation of an almost generic curve
is called \emph{fiberwise} if the curve remains
almost generic throughout the deformation.

\emph{Deformation equivalence} of real trigonal curves is
the equivalence relation generated by equivariant
fiberwise deformations and
real isomorphisms (in the category of pairs as above).

\subsection{Auxiliary statements}\label{auxiliary}
Let $\pi\:\Sigma\to B$ and $E\subset\Sigma$ be as in
Subsection~\ref{real-trigonal}.
Recall that, for any
coefficient group~$G$, the inverse
Hopf homomorphism $\pi^*$ establishes an isomorphism
$$
\pi^*\:H_1(E;G)@>\cong>>H_3(\Sigma;G).
\eqtag\label{eq.tr}
$$

Let~$C \subset \Sigma$ be a nonsingular trigonal curve.
Assume that $d=2k$
is even and
consider a double covering $p\:X\to\Sigma$ of~$\Sigma$ ramified
at $C+E$. It is a Jacobian elliptic surface. Let
$\omega\in H^1(\Sigma\sminus(C\cup E);\Z_2)$ be the class of the
covering and denote by $\tr$ the transfer homomorphism
$$
\tr\:H_*(\Sigma,C\cup E;\Z_2)\to H_*(X;\Z_2).
$$

\lemma\label{tr}
The composition
$$
\def\ZZZ{}
H_1(E\ZZZ)@>\pi^*>>H_3(\Sigma\ZZZ)
@>\rel>>H_3(\Sigma,C\cup E\ZZZ)
@>{}\cap\omega>>H_2(\Sigma,C\cup E\ZZZ)
@>\partial>>H_1(C\ZZZ)\oplus H_1(E\ZZZ)
$$
\rom(all homology with coefficients $\Z_2$\rom) is given by
$a\mapsto\pi_C^*a\oplus a$, where
$\pi_C^*$ stands for the inverse
Hopf homomorphism $H_1(E;\Z_2)\to H_1(C;\Z_2)$.
\endlemma

\proof
Realize a class in $H_1(B;\Z_2)$ by an embedded circle
$\gamma\subset B$ and restrict all maps to~$\gamma$ to obtain
$X_\gamma\to\Sigma_\gamma\to\gamma$. Then
$\pi^*[\gamma]=[X_\gamma]$ and $\rel[X_\gamma]\cap\omega$ is the
class dual to~$\omega$; its boundary is the fundamental class of
the ramification locus.
\endproof

\corollary\label{dKer}
One has
$$
\def\ZZZ{}
\partial\Ker[\tr\:H_2(\Sigma,C\cup E\ZZZ)\to H_2(X\ZZZ)]=
 \Im[\pi_C^*\oplus\id\:H_1(E\ZZZ)\to H_1(C\ZZZ)\oplus H_1(E\ZZZ)]
$$
\rom(all homology with coefficients $\Z_2$\rom).
\endcorollary

\proof
Comparing the Smith exact sequence
$$
\CD
H_3(\Sigma,C\cup E)@>\omega\oplus\partial>>
H_2(\Sigma,C\cup E)\oplus H_2(C\cup E)@>\tr+\tilde{\inj}_*>>
H_2(X)
\endCD
$$
of the double covering~$p$ (where all homology groups are with
coefficients~$\Z_2$ and
$\tilde{\inj}_*\:H_2(C\cup E;\Z_2)\to H_2(X;\Z_2)$ is the inclusion
homomorphism) and the exact sequence of the pair
$(\Sigma,C\cup E)$, one concludes that $\Ker\tr$ is the image of
the composed homomorphism
$$
\def\ZZZ{;\Z_2}
H_3(\Sigma\ZZZ)
@>\rel>>H_3(\Sigma,C\cup E\ZZZ)
@>{}\cap\omega>>H_2(\Sigma,C\cup E\ZZZ).
$$
Hence, the statement follows from the isomorphism~\eqref{eq.tr}
and Lemma~\ref{tr}.
\endproof

The following statement is well known.

\lemma\label{w2}
For a Jacobian elliptic surface $p\:X\to\Sigma$ ramified at $C+E$
one has $w_2(X)=kp^*(f)\in H_2(X;\Z_2)$ and
$p^*(e)=0\in H_2(X;\Z_2)$.
\endlemma

\proof
Let~$\tC$ and~$\tE$ be the pull-backs of~$C$ and~$E$,
respectively, in~$X$. Since the group $H_1(X)=H_1(B)$
is torsion free, so is $H_2(X)$ and one
has
$$
[\tC]+[\tE]=\frac12p^*([C]+[E])=p^*(2e+3kf).
$$
(Recall that, for an algebraic
curve $D\subset\Sigma$, one has $p^*[D]=[p^*D]$, where $p^*D$
is the divisorial pull-back of~$D$. The reduction $p^*\bmod2$ is the
composition of the relativization
$\rel\:H_2(\Sigma;\Z_2)\to H_2(\Sigma,C\cup E;\Z_2)$
and the transfer~$\tr$.)
Then, due to the projection formula, one has
$$
w_2(X)=p^*w_2(\Sigma)+[\tC]+[\tE]=kp^*(f)
$$
in $H_2(X;\Z_2)$.
For the last statement, $p^*(e)=2[\tE]=0\bmod2$.
\endproof

\subsection{Trigonal curves of type~$\Ib$}\label{s.Ib}
Recall that a nonsingular real curve~$C$ with nonempty real part
is said to be of \emph{type~$\I$}, or \emph{separating}, if
$[C_\R] = 0 \in H_1(C; \Z_2)$;
otherwise, $C$ is of \emph{type~$\II$}.

If~$C$ is a (connected) separating real curve, the complement $C
\sminus C_\R$ splits into two connected components. Their closures
are called \emph{halves} of~$C$ and denoted~$C^\pm$.
One has $C_\R =
\partial C^+ = \partial C^-$.

In the case of real trigonal curves in a real ruled surface
$\pi\:\Sigma\to B$, one can consider a wider class that shares
most useful properties of curves of type~$\I$.

\definition\label{typeIB}
A real trigonal curve $C\subset\Sigma$ is said to be of
\emph{type~$\Ib$} if the identity
$[\CR]=\pi_C^*[\BR]$
holds in $H_1(C;\Z_2)$.
\enddefinition

\lemma\label{base.type.I}
A trigonal curve $C\subset\Sigma$ is of type~$\I$ if and only if
$C$ is of type~$\Ib$ and the base~$B$ is of type~$\I$.
\endlemma

\proof Clearly, types~$\I$ and~$\Ib$ are equivalent whenever $B$ is
of type~$\I$. Hence, it suffices to prove that the base~$B$ of a
trigonal curve~$C$ of type~$\I$ is necessarily of type~$\I$.
Represent~$C$ as the union of two halves, $C=C^+\cup C^-$, and
define functions $n^\pm\:B\to\Z$ \via\
$n^\pm(b)=\#(\pi_C^{-1}(b)\cap C^\pm)$.
On the complement $B^\circ\sminus\BR$ both functions
$n^\pm$ are locally constant and one has
$n^++n^- = 3$
and
$c^*n^\pm=n^\mp$ (since $C^+$ and~$C^-$ are interchanged by the real
structure on~$C$). Hence, one can define a half~$B^+$ of~$B$ as the
closure of the set $\{b\in B\,|\,n^+(b)<n^-(b)\}$.
\endproof

Consider a trigonal curve~$C$ of type~$\Ib$ and define $\CIm$ as
the closure of the set $\pi_C^{-1}(\BR)\sminus\CR$. Let
$\BIm=\pi_C(\CIm)$. Clearly, $\CIm=\varnothing$ if and only if $C$
is hyperbolic.
The following statement is immediate.

\lemma\label{Type.Ib}
A real trigonal curve~$C$ is of type~$\Ib$ if and only if the
class $[\CIm]$ vanishes in $H_1(C;\Z_2)$.
\qed
\endlemma

\paragraph\label{ss.C+}
According to the previous lemma, a non-hyperbolic trigonal curve~$C$
of type~$\Ib$ can be represented as the union of two surfaces~$C_+$
and~$C_-$ (possibly disconnected), disjoint except for the common
boundary $\partial C_+=\partial C_-=\CIm$. Define functions
$m_\pm\:B\to\Z$ \via\ $m_\pm(b)=\#(\pi_C^{-1}(b)\cap
C_\pm)-\chi\subIm(b)$, where
$\chi\subIm$ is the characteristic function of $\BIm$. It is easy to
see that, on the subset $B^\circ\subset B$, both functions~$m_\pm$
are locally constant and one has $m_++m_-=3$. Since $B^\circ$ is
connected, $m_\pm|_{B^\circ}=\const$. In what follows, we mark the
surfaces~$C_\pm$ so that $m_+|_{B^\circ}\equiv1$ and
$m_-|_{B^\circ}\equiv2$.

Due to the convention above, the restriction $\pi_+\:C_+\to B$
of~$\pi_C$ is one-to-one except on the boundary~$\partial C_+$. In
particular, it follows that $C_+$ is connected unless $B$ is of
type~$\I$ and $\BR=\BIm$. In any case, both~$C_+$ and~$C_-$ are
invariant under the real structure on~$C$.

\subsection{Jacobian surfaces}\label{s.If}
A real surface~$X$ is said to be of \emph{type~$\I$} if
$[\XR]=w_2(X)$ in $H_2(X;\Z_2)$.
A real elliptic surface~$X$ is said
to be of \emph{type~$\If$} if the image of~$[\XR]$ in
$H_2(X;\Z_2)$ is a multiple of the class of a fiber of~$X$ (\cf.
Lemma~\ref{w2}).
Fix a real ruled surface $\pi\:\Sigma\to B$ over a
real base~$B$ and assume
that the self-intersection of the exceptional section
$E\subset\Sigma$ is even, $E^2=-2k$.

\lemma\label{IF=I}
A Jacobian real elliptic surface~$X$ is of type~$\I$ if and only if it
is of type~$\If$.
\endlemma

\proof Let $\tE$ be the real section of~$X$.
Then
$[\XR]\circ[\tE]=k=kp^*(f)\circ[\tE]$, and it remains to apply
Lemma~\ref{w2}.
\endproof

\proposition\label{IF=IB}
Let $C\subset\Sigma$ be a real trigonal curve. Then, a
Jacobian elliptic surface~$X$ ramified at $C+E$ is of
type~$\I$ if and only if $C$ is of type~$\Ib$.
\endproposition

\proof
In view of Lemma~\ref{IF=I}, it suffices to show that $X$ is of
type~$\If$ if and only if $C$ is of type~$\I$.
Recall that the class $[\XR]\in H_2(X;\Z_2)$
can be represented in the form
$\tr[\SigmaR^+]$, where $\SigmaR^+\subset\SigmaR$ is the appropriate
half of the real part $\SigmaR$
and $[\SigmaR^+]$ is regarded as a relative
class in $H_2(X,C\cup E;\Z_2)$.
Then,
the identity $[\XR]=ap^*(f)$ holds for some $a\in\Z_2$ if and only
if $\tr([\SigmaR^+]-af)=0$. Since $\partial f=0$ and $p^*(e)=0$, see
Lemma~\ref{w2}, the statement follows from
Corollary~\ref{dKer}.
\endproof


\subsection{Other surfaces}
Recall that, to every elliptic surface~$X$, one can assign its
\emph{Jacobian surface}~$J$. If $X$ is real, then $J$ also
inherits a canonical real structure.

%

\corollary[Conjecture]
A real elliptic surface~$X$ is of type~$\If$
if and only if the real trigonal curve constituting the ramification
locus of the Jacobian surface of~$X$ is of type~$\Ib$.
\endcorollary

\section{Dessins}\label{rtc}

The notion of dessin used in this paper is a real version of
Grothendieck's \emph{dessins d'enfants}, adjusted for the study of
real meromorphic functions with a certain preset ramification over
the three real points $0,1,\infty\in\Cp1$. More precisely, we
consider the quotient by the complex conjugation of a properly
decorated pull-back of $(\Rp1;0,1,\infty)$, the pull-backs of~$0$,
$1$, and~$\infty$ being marked with \black--, \white--, and
\cross--, respectively. Note that, unlike Grothendieck's original
setting, the functions considered
may (and usually do) have other critical values, which
are ignored unless they are real.

In the exposition below we follow~\cite{DIK},
omitting most proofs and references.

\subsection{Trichotomic graphs}\label{tg}
Let $\base$ be a (topological) compact connected surface, possibly
with boundary.
(In the topological
part of this section we are working in the {\sl PL\/}-category.)
We use the term \emph{real} for points, segments, \etc.
situated at the boundary~$\partial\base$. For a graph
$\Gamma\subset\base$, we denote by~$\Base\Gamma$ the closed cut
of~$D$ along~$\Gamma$.
The connected components of~$\Base\Gamma$
are called \emph{regions} of~$\Gamma$.

A \emph{trichotomic graph} on~$\base$
is an embedded oriented
graph $\Gamma\subset\base$
decorated with the
following additional structures (referred to as \emph{colorings}
of the edges and vertices of~$\Gamma$, respectively):
\roster
\item"--"
each edge of~$\Gamma$ is of one of the three kinds: \solid-, \bold-,
or \dotted-;
\item"--"
each vertex of~$\Gamma$ is of one of the four kinds: \black-, \white-,
\cross-, or monochrome (the vertices of the first three kinds being called
\emph{essential});
\endroster
and satisfying the following conditions:
\roster
\item\local{tg-boundary}
the boundary $\partial\base$ is a union of
edges
and vertices
of~$\Gamma$;
\item\local{tg-valency}
the valency of each essential vertex of~$\Gamma$ is at
least~$2$,
and the valency of each monochrome vertex of~$\Gamma$ is at
least~$3$;
\item\local{tg-oriented}
the orientations of the edges of~$\Gamma$ form an orientation
of the boundary $\partial\Base\Gamma$;
this orientation extends to an orientation of~$\Base\Gamma$;
\item\local{tg-monochrome}
all edges incident to a monochrome vertex are of the same kind;
\item\local{tg-cross}
\cross-vertices are incident to incoming \dotted-edges and
outgoing \solid-edges;
\item\local{tg-black}
\black-vertices are incident to incoming \solid-edges and
outgoing \bold-edges;
\item\local{tg-white}
\white-vertices are incident to incoming \bold-edges and
outgoing \dotted-edges;
\item\local{tg-triangle}
each \emph{triangle} (\ie, region with three essential vertices in
the boundary) is a topological disk.
\endroster
In \loccit{tg-cross}--\loccit{tg-white} the lists are complete,
\ie, vertices cannot be incident to edges of other kinds or with
different orientation.

In view of~\loccit{tg-monochrome}, the monochrome vertices can
further be subdivided into \solid-, \bold-, and \dotted-,
according to their
incident edges.
A \emph{monochrome cycle} in~$\Gamma$ is a cycle with all vertices
monochrome, hence all edges and vertices of the same kind.


\paragraph\label{full-valency}
Let~$B$ be the oriented double of~$\base$,
and denote by $\Gamma' \subset B$
the
preimage of~$\Gamma$,
with each vertex and edge decorated according to its image
in~$\Gamma$.
(Note that $\Gamma'$ is also a trichotomic graph.)
The \emph{full valency} of a vertex of~$\Gamma$ is
the valency of any preimage of this vertex
in~$\Gamma'$.
The full valency of an inner
vertex coincides with its valency in~$\Gamma$; the full valency of a real
vertex equals $2\cdot\text{valency}-2$.
Conditions~\loccit{tg-oriented} and~\loccit{tg-boundary}
imply
that the
orientations of the edges of~$\Gamma'$ incident to a vertex alternate.
Thus,
the full valency of any
vertex is even.

\paragraph\label{pillar}
The collection of all vertices and edges of
a trichotomic graph~$\Gamma$
contained in a given connected component
of~$\partial\base$ is called a {\it real component}
of~$\Gamma$. In the drawings, (portions of) the real components
are indicated by wide grey lines.
A real component
(and the corresponding
component
of $\partial\base$) is
called
\Dashes
\dash \emph{even}/\emph{odd}, if it contains
an even/odd number of \white-vertices of~$\Gamma$,
\dash \emph{hyperbolic}, if all
edges of this component are dotted,
\dash \emph{anti-hyperbolic}, if the component
contains no dotted edges.
\endDashes
A trichotomic graph
is called \emph{hyperbolic} if all its real components are
hyperbolic.

If a union of (the closures of)
some real edges of the same kind is homeomorphic to
a closed interval,
this union is called a \emph{segment}.
A \dotted- (\bold-)
segment
is called \emph{maximal}
if it is bounded by
two
\cross--
(respectively, \black--) vertices.
Define the \emph{parity} of a maximal
segment as the parity
of the number of \white-vertices contained in the segment.
A \emph{pillar} is
either a hyperbolic component,
or a maximal \dotted- or \bold-segment.

\subsection{Dessins}\label{sketches}
Recall that to any trigonal curve~$C \subset \Sigma$,
$\Sigma \to B$, one can associate its
\emph{\rom(functional\rom) $j$-invariant} $j = j_C \:B\to\Cp1$,
which is the analytic continuation
of the
meromorphic function $B^\circ\to\C$ sending each
fiber $F \subset \Sigma$ nonsingular for $C$
to the conventional $j$-invariant (symmetrized cross-ratio)
of the quadruple of points $(C \cup E) \cap F$
in the
projective line~$F$; following Kodaira, we divide the
$j$-invariant by $12^3$, so that its `special' values are $j=0$
and~$1$ (corresponding to quadruples
with a symmetry of order~$3$ and~$2$, respectively).

We assume that the target Riemann sphere $\Cp1=\C\cup\{\infty\}$
is equipped with the standard real structure $z \to \bar z$.
With respect to this real structure, the $j$-invariant of
a real trigonal curve is real, descending to a map
from the quotient $D = B/c$ to the disk $\Cp1\!/\spbar{}$.
The pull-back of the real part $\Rp1=\partial(\Cp1\!/\spbar)$
under this map
is denoted by~$\Gamma_C$.
This pull-back, regarded as a graph in~$D$,
has a natural trichotomic
graph structure:
the \black--, \white--, and
\cross-vertices are  the pull-backs of~$0$, $1$,
and~$\infty$, respectively (monochrome vertices being the
branch points with other real critical values), the edges
are \solid-, \bold-, or \dotted- provided that their images belong to
$[\infty,0]$, $[0,1]$, or $[1,\infty]$, respectively, and
the
orientation of~$\Gamma_C$ is that induced from the positive
orientation of~$\Rp1$ (\ie, order of~$\R$).
This definition implies
that~$\Gamma_C$ has no oriented monochrome cycles.
Furthermore, the boundary of each triangle is mapped to $\Rp1$
with degree one, and any extension of this map to the triangle
itself also has degree one; hence, the triangle is homeomorphic to
$\Cp1\!/\spbar{}$, which explains
condition~\itemref{tg}{tg-triangle}.

A real
trigonal curve~$C$ is almost generic if and only if
the full valency of each \cross--, \white--, and
\black-vertex
of~$\Gamma_C$ equals, respectively,~$2$, $0\bmod4$, and $0\bmod6$.
A real trigonal curve~$C$ is called \emph{generic} if
its graph $\Gamma = \Gamma_C$
has the following properties:
\roster
\item\local1
the full valency of each \cross--, \white--, or
\black-- vertex of~$\Gamma$ is, respectively,~$2$, $4$, or $6$;
\item\local2
the valency of any real monochrome vertex of~$\Gamma$ is~$3$;
\item\local3
$\Gamma$ has no inner monochrome vertices.
\endroster
A trichotomic graph
satisfying conditions~\loccit1--\loccit3
and without oriented monochrome cycles
is called a \emph{dessin}.
We freely extend to dessins
all
terminology that
applies to
generic trigonal curves.
Since we only consider curves with nonempty real part,
we always assume that the boundary of the
underlying surface of a dessin is nonempty.

Any almost generic
real
trigonal curve can be perturbed to a generic one.



\paragraph\label{ss.equivalence}
Two dessins are called
\emph{equivalent} if, after a homeomorphism
of the underlying surfaces, they
are
connected by a finite sequence of
isotopies and the following \emph{elementary moves}:
\roster
\item"--"
\emph{monochrome modification}, see
Figure~\ref{fig.moves}(a);
\item"--"
\emph{creating \rom(destroying\rom) a bridge}, see
Figure~\ref{fig.moves}(b),
where a \emph{bridge} is a pair of
monochrome vertices connected by a real monochrome edge;
\item"--"
\emph{\white-in} and its inverse \emph{\white-out}, see
Figure~\ref{fig.moves}(c) and~(d);
\item"--"
\emph{\black-in} and its inverse \emph{\black-out}, see
Figure~\ref{fig.moves}(e) and~(f).
\endroster
(In the first two cases, a move is
considered valid
only if the result
is again a dessin.
In other words,
one needs to check
the absence of oriented monochrome cycles
and triangular regions other than disks.)
An equivalence of two dessins in the same underlying surface~$\base$
is called \emph{restricted} if the homeomorphism
is identical and the isotopies above can be chosen
to preserve the pillars (as sets).

\midinsert
\eps{Fig1all}
\figure\label{fig.moves}
Elementary moves of dessins
\endfigure
\endinsert

\Remark\label{rem.orientation} In view of
Condition~\iref{tg}{tg-oriented} in the definition of trichotomic
graph, any monochrome modification and
creation/destruction of a bridge automatically
respect the orientations of the edges involved, see
Figure~\ref{fig.moves}.
This fact is in a
contrast with the definition of equivalence of skeletons, see
Subsection~\ref{equivalence} below, where respecting a certain
orientation is an extra requirement.
\endRemark

The following statement is proved in~\cite{DIK}.

\proposition\label{equiv.curves}
Each dessin~$\Gamma$ is of the form~$\Gamma_C$
for some generic real trigonal curve~$C$.
Two generic real trigonal curves are deformation
equivalent
\rom(in the
class of
almost generic
real trigonal curves\rom)
if and only if their dessins are equivalent.
\qed
\endproposition

\paragraph\label{tg->topology}
The definition of the $j$-invariant gives one an easy way to
reconstruct the topology of a generic real trigonal curve
$C\subset\Sigma$ from its dessin~$\Gamma_C$. Let
$\pi\:\Sigma\to B$ and $\pi_C$ be as in
Subsection~\ref{real-trigonal}.
Topologically,
the base~$B$ is the orientable double of the underlying
surface~$\base$ of~$\Gamma_C$.
Let $\Gamma' \subset B$
be the decorated
preimage of~$\Gamma_C$, see~\ref{full-valency}.
Then
$B^\circ = B \sminus \{\text{\cross-vertices of $\Gamma'$}\}$
and the pull-back
$\pi_C^{-1}(b)$ of a point
$b\in B^\circ$ consists of three
points in the complex affine line $\pi^{-1}(b) \sminus E$.
\roster
\item\local{b.region}
If $b$ is an inner point of a region of~$\Gamma'$, the three points
of $\pi_C^{-1}(b)$ form a triangle~$\Delta_b$ with all three
edges distinct. As a consequence, the restriction of~$\pi_C$ to the
interior of each region of~$\Gamma'$ is a trivial covering.
\item\local{b.dotted}
If $b$ belongs to a \dotted- edge of~$\Gamma'$, the three points of
$\pi^{-1}_C(b)$ are collinear. The ratio
$(\text{smallest distance})/(\text{largest distance})$
is in $(0,1/2)$; it
tends to~$0$ (respectively, $1/2$) when $b$ approaches a \cross--
(respectively,~\white--) vertex.
\item\local{b.edge}
If $b$ belongs to a \solid- (\bold-) edge of~$\Gamma'$, the three
points of $\pi^{-1}_C(b)$ form an isosceles triangle
with the angle at the vertex less than (respectively, greater than)
$\pi/3$. The angle
tends to~$0$, $\pi/3$, or~$\pi$ when $b$ approaches, respectively,
a \cross--, \black--, or \white-vertex.
\endroster

The number of \white-vertices of~$\Gamma'$
is called the \emph{degree} $\deg\Gamma$ of~$\Gamma$.
One has $\deg\Gamma = 0 \bmod 3$,
and $-\frac{1}{3}\deg\Gamma = E^2$,
where~$E \subset \Sigma$ is the exceptional section.

\paragraph\label{ss.labels}
In view of~\iref{tg->topology}{b.region},
over the interior of each region~$R$ of
the pull-back $\Gamma'\subset B$
there is a canonical way to label the three sheets of~$C$ by~\1, \2,
\3, according to the increasing of the opposite side of
the
triangle~$\Delta_b$
over any point $b\in R$.
This labelling is obviously preserved by the real structure $c\:B\to
B$ and hence descends to the regions of~$\Gamma$. The passage
through an edge of~$\Gamma$ results in the following transformation:
\Dashes \dash solid edge: the transposition $(23)$; \dash bold edge: the transposition
$(12)$; \dash dotted edge: the change of the orientation
of~$\Delta$.
\endDashes
The transpositions above represent a change of the labelling
rather than a nontrivial monodromy. Although the change of
orientation of~$\Delta$ makes sense, its orientation itself is
only well defined if $B$ is of type~$\I$ and a half of~$B$ is
chosen.

\paragraph\label{tg->real}
The real components $\Gamma_i$ of~$\Gamma$
are identified with the connected components~$B_i$
of~$B_\R$.
The pull-back $\pi_C^{-1}(b)$ of a real point $b\in\partial\base$
has three real points if $b$ is a \dotted- point or a
\white-vertex
adjacent to two real \dotted- edges; it has two real points, if $b$ is a
\cross-vertex, and a single real point otherwise.
A component $\Sigma_i$ of $\Sigma_\R$ is orientable
if and only if the corresponding real
component $\Gamma_i$ is even.

A component~$B_i$ is (anti-)hyperbolic
(see~\ref{real-trigonal})
if and only if so is~$\Gamma_i$.
If $B_i$ is non-hyperbolic,
its ovals and zigzags are represented
by the maximal \dotted- segments
of~$\Gamma_i$, even and odd, respectively.
The latter are also called \emph{ovals} and \emph{zigzags}
of~$\Gamma$.

\section{Trigonal curves of type~$\Ib$}\label{curves.typeIb}

In this section, we characterize the dessins of trigonal curves
of type~$\Ib$ and study their basic properties.

\subsection{Canonical labelling}\label{canonical-label}
Let $C$ be a \emph{non-hyperbolic} trigonal curve of type~$\Ib$,
and let $C_+\subset C$ be the surface mapped generically
one-to-one to~$B$, see~\ref{ss.C+}.


Let $\Gamma\subset\base$ be the dessin of~$C$.
Since $C_+$ is $c$-invariant,
each region~$R$ of~$\Gamma$ can be labelled according
to the label of
the sheet of~$C_+$
over~$R$, see~\ref{ss.labels}.
Then,
each \emph{inner}
edge~$e$ of~$\Gamma$ can be labelled according to the label(s) of
the adjacent regions.
The possible labels are as follows:
\Dashes
\dash
an inner solid edge can be of type~\1 or~\n1 (not~\1);
\dash
an inner bold edge can be of type~\3 or~\n3 (not~\3);
\dash
an inner dotted edge can be of type~\1, \2, or~\3.
\endDashes
(One cannot distinguish types~\2 and~\3 along a solid edge or
types~\1 and~\2 along a bold edge due to the relabelling mentioned
above;
in these cases, we assign to the edges types~\n1 and~\n3,
respectively.)
On the contrary,
the same rule assigns a
well defined label~\1, \2, or~\3 to
each \emph{real}
edge~$e$ of~$\Gamma$: the relabelling
in~~\ref{ss.labels}
is compensated for by the discontinuity of~$C_+$
across~$\BIm$.

\lemma\label{no13}
A real solid edge cannot be of type~\1\rom; a real bold edge
cannot be of type~\3.
\endlemma

\proof
Otherwise, the surface $C_+\to B$ would be two-sheeted over the
regions of~$\Gamma$ adjacent to the edge.
\endproof

\theorem\label{type.distrib} A non-hyperbolic generic trigonal
curve~$C$ is of type~$\Ib$ if and only if the regions of its dessin
$\Gamma \subset \base$ admit a labelling which
satisfies the following
conditions:
\roster
\item\local{regions1}
the region adjacent
to a real solid \rom(bold\rom) edge is
not of type~$\1$ \rom(respectively, not of type~$\3$);
\item\local{regions2}
the two regions adjacent
to an inner solid edge are either both of type~$\1$ or
of distinct types~$\2$ and~$3$;
\item\local{regions3}
the two regions adjacent
to an inner bold edge are either both of type~$\3$ or
of distinct types~$\1$ and~$2$;
\item\local{regions4}
the two regions adjacent
to an inner dotted edge
are of the same type.
\endroster
\endtheorem

\proof
If~$C$ is of type~$\Ib$, its labelling defined above
does satisfy~\loccit{regions1}--\loccit{regions4}:
Property \loccit{regions1} is the statement of Lemma~\ref{no13},
and Properties~\loccit{regions2}--\loccit{regions4}
follow from~\ref{ss.labels}.

For the converse,
lift the labelling to the preimage $\Gamma' \subset B$
of~$\Gamma$, \cf.~\ref{tg->topology}, over each region of
$\Gamma'$ take the sheet selected by the labelling, and
define~$C_+$ as the closure of the union of these sheets.
Then, in view of~\loccit{regions1}--\loccit{regions4},
one has $\partial C_+=\CIm$, \ie,
$[\CIm] = 0 \in H_1(C;\Z_2)$, and $C$ is of type~$\Ib$ due to
Lemma \ref{Type.Ib}.
\endproof

\subsection{Dessins of type~$\I$}\label{dessins-typeI}
A dessin~$\Gamma$ equipped with a labelling satisfying
Conditions~\iref{type.distrib}{regions1}--\ditto{regions4}
is said to be \emph{of type~$\I$}. We assume the labelling
extended to edges as explained in Subsection~\ref{canonical-label}.
Fix a dessin~$\Gamma$ of type~$\I$.
Below, we discuss further properties of its labelling
and extend it to some other objects related to~$\Gamma$.

\lemma\label{inner.vertex}
The edges
adjacent to an inner vertex
are labelled as follows\rom:
\Dashes
\dash
\cross-vertex\rom: $(\1,\1)$\rom;
\dash
\black-vertex\rom: $(\1,\n3,\n1,\3,\n1,\n3)$\rom;
\dash
\white-vertex\rom: $(\3,\3,\3,\3)$ or $(\n3,\1,\n3,\2)$.
\qed
\endDashes
\par\removelastskip
\endlemma

\lemma\label{real.vertex}
The edges adjacent to a real vertex are labelled as follows\rom:
\Dashes
\dash
\cross-vertex\rom: both edges are of the same type~\2 or~\3\rom;
\dash
\black-vertex\rom: $(\2,\n3,\1,\1)$ or $(\3,\3,\n1,\2)$\rom;
\dash
\white-vertex with real edges dotted\rom:
$(\3,\3,\3)$ or $(\1,\n3,\2)$\rom;
\dash
\white-vertex with real edges bold\rom: all edges are of the
same type~\1 or~\2.
\qed
\endDashes
\par\removelastskip
\endlemma

According to Lemmas~\ref{inner.vertex} and~\ref{real.vertex}, the
two edges adjacent to a single \cross-vertex are always of the
same type. We assign this type to the vertex itself, thus speaking
about \emph{\cross-vertices of type~\1} (necessarily inner)
\emph{or~\2 or~\3} (necessarily real).

\corollary\label{oval.type}
The two \cross-vertices bounding a single oval of~$\Gamma$ are of
the same type, which can be~\3 or~\2.
In the former case, all dotted edges constituting the oval are
of type~\3\rom; in the latter case, the type of the edges
alternates between~\2 and~\1 at each \white-vertex.
The two \cross-vertices
bounding a single zigzag of~$\Gamma$ are always of type~\3, and so
are all dotted edges constituting a zigzag.
\qed
\endcorollary

According to the types of the dotted edges constituting an
oval, we will speak about ovals of type~\3 and~\2 (if
all dotted edges are of
type~\3 or~\2, respectively) and ovals of type~\n3 (if there are
edges both of type~\2 and~\1). Note that ovals of type~\2 are
necessarily `short', \ie, they contain no \white-vertices, whereas
each oval of type~\n3 necessarily contains a \white-vertex.

By definition, each zigzag is regarded to be of type~\3.

\lemma\label{dotted.type}
The real dotted edges
constituting an odd
\rom(respectively, even\rom)
hyperbolic component of~$\Gamma$
are all of type~\3 \rom(respectively, are all of the same
type~\3, \2, \1 or
alternate between type~\2 and~\1
at each \white-vertex\rom).
\endlemma

\proof
Within each real dotted segment, the types of the edges either are
$\const=\3$ or alternate between~\2 and~\1 at each \white-vertex.
In the latter case, the number of vertices must be even.
\endproof

According to the types of the dotted edges constituting the
component, we will speak about
hyperbolic components of type~\3, \2,
\1, or~\n3. Note that components of type~\2 or~\1
do not contain \white-vertices, whereas each component
type~\n3 necessarily contains a \white-vertex.

\lemma\label{no.bold}
The real part
of~$\Gamma$ has no odd anti-hyperbolic components.
\endlemma

\proof After destroying solid bridges and
a sequence of \black-ins along inner solid edges, one
can assume that the edges constituting an anti-hyperbolic component
are all bold. They cannot be of type~\3, see Lemma~\ref{no13};
hence, their types alternate between~\1 and~\2 at each inner bold
edge attached to the component, and the component must be even.
\endproof

\subsection{Unramified dessins of type~$\I$}
A dessin is called
\emph{unramified},
if all its \cross-vertices
are real.
In other words, unramified are the dessins corresponding
to maximally inflected curves.
In this subsection, we assume that~$\Gamma$
is an unramified dessin of type~$\I$.

\lemma\label{no1}
The dessin~$\Gamma$ has no solid or dotted edges of type~\1.
\endlemma

\proof
A solid or dotted edge of type~\1 would end at a \cross-vertex
of type~\1
(possibly, passing through a number of monochrome vertices), which
would have to be inner.
\endproof

%

\lemma\label{black.real}
Each \black-vertex $v\in\Gamma$ is real, and the edges adjacent
to~$v$ are of types $(\3,\3,\n1,\2)$. The immediate essential neighbors
of~$v$ in the real part of~$\Gamma$ are a \white-- and a
\cross-vertex.
\endlemma

\proof
The first two statements follow from Lemmas~\ref{inner.vertex},
\ref{real.vertex}, and~\ref{no1}. If $v$ had another
\black-vertex as an immediate essential neighbor,
the two vertices could be
pulled in by a \black-in transformation, producing an inner
\black-vertex.
\endproof

\corollary\label{no.not3}
The dessin~$\Gamma$ has no inner bold edges of type~\n3.
\qed
\endcorollary

\corollary\label{white.monotype}
All edges adjacent to a single \white-vertex of~$\Gamma$ are of
the same type. \rom(Accordingly, we will speak about the
\emph{type of a \white-vertex}.\rom) A real \white-vertex with
real edges bold is of type~\2\rom; all other \white-vertices
are of type~\3.
\endcorollary

\proof
The types of edges adjacent to a \white-vertex are listed in
Lemmas~\ref{inner.vertex} and \ref{real.vertex}, and all but a few
possibilities are eliminated by Lemma~\ref{no1}.
\endproof

\corollary\label{short.oval}
The dessin~$\Gamma$ has no ovals of type~\n3.
\qed
\endcorollary

\lemma\label{white.real}
Let $v\in\Gamma$ be a \white-vertex of type~\2.
Then $v$ is real, and the immediate essential neighbors of~$v$
in the real part of~$\Gamma$
are two \black-vertices.
\endlemma

\proof
According to Corollary~\ref{white.monotype},
the vertex~$v$ is real
and the real edges adjacent to~$v$ are bold.
Hence, the neighbors of~$v$ are either
\white-- or \black-vertices.
If another \white-vertex~$u$ (possibly, $v$ itself)
were a neighbor of~$v$, the
dotted
segment
connecting~$u$ and~$v$ would contain a monochrome vertex with
an inner bold edge of type~\n3 adjacent to it. This would contradict to
Corollary~\ref{no.not3}.
\endproof

A pillar
consisting of a
\white-vertex and pair of real bold segments connecting it to
\black-vertices, as in Lemma~\ref{white.real}, is called a
\emph{jump}. To each jump, we assign type~\2, according to the
types of its \white-vertex and bold edges.

\proposition\label{no.monochrome}
Any hyperbolic component of~$\Gamma$ is of type~\3 or~\2 \rom(and
in the latter case it is free of \white-vertices\rom). Any
anti-hyperbolic component of~$\Gamma$
is formed by solid edges and solid
monochrome vertices.
\endproposition

\proof
As in the proof of Lemma~\ref{no.bold},
if an anti-hyperbolic component has a bold edge,
one can assume all edges of this component bold.
Any hyperbolic
component of type~\n3 or any real component
with all edges bold
would have an inner bold edge of type~\n3
attached to it; this contradicts to Corollary~\ref{no.not3}.
Similarly, any
hyperbolic
component of type~\1 would have an inner
dotted edge of type~\1 attached to it; this contradicts to
Lemma~\ref{no1}.
\endproof

The following theorem summarizes the results of this section.

\theorem\label{th.summary}
Let $\Gamma$ be an unramified dessin of type~$\I$. Then
\roster
\item
the pillars of~$\Gamma$ are ovals, zigzags, jumps,
and hyperbolic components\rom;
\item
each pillar has a well defined \emph{type}, \2
or~\3, all jumps being of type~\2 and all zigzags being of
type~\3\rom;
\item\local1
pillars of type~\2
are interconnected by inner dotted edges of type~\2\rom;
these edges, as well as
pillars
of type~\2 other than jumps, are
free of \white-vertices\rom;
\item\local2
pillars of type~\3
are interconnected by inner dotted edges of type~\3
or pairs of such edges attached
to an inner \white-vertex each\rom;
\item\local3
the following \emph{parity rule} holds\rom:
along each real component of~$\Gamma$,
the types
of the pillars alternate.
\qed
\endroster
\par\removelastskip
\endtheorem

\subsection{Complex orientations}\label{c.orientaions}
Recall that the real part $C_\R$ of any connected
algebraic curve~$C$ of type~$\I$
admits a distinguished pair of opposite orientations,
called \emph{complex orientations}, which are induced
on the common boundary $C_\R=\partial C^\pm$
by the
complex orientations of the two halves $C^\pm$ of~$C$.

Let $C \subset \Sigma$ be a nonsingular
real trigonal curve
over a base~$B$ of type~$\I$. Consider the set $\Cal B \subset B_\R$
of real fibers of~$\Sigma$ that intersect $C_\R$ in a single
point each (counted with multiplicity). Denote by $\bar\Cal B$
the closure of~${\Cal B}$, and let~$\bar L$ be the restriction
to $\bar{\Cal B}$ of the real ruling
$\pi\:\Sigma_R \sminus E_\R \to B_\R$. It is a real affine
line bundle.
 Any
orientation of~$C_\R$ induces in an obvious way an orientation of
the restriction $\bar L |_{\bar \Cal B \sminus \Cal B}$.

\definition\label{def.qc.orientations}
For a non-hyperbolic curve~$C$, an orientation of~$C_\R$ is called
\emph{quasi-complex} if the induced orientation of $\bar L |_{\bar
\Cal B \sminus \Cal B}$ extends to~$\bar L$ and, with respect to
some complex orientation of~$B_\R$, the restriction of the
projection $\pi_C\:C_\R \to B_\R$ is of degree~$+1$ over each
component of~$B_\R$.
\enddefinition

\Remark\label{remark.qc.orientations} Denote by~$Z$ the projection
to $B_\R$ of the union of all zigzags and real vertical flexes
of~$C_\R$. Then, $C_\R$ admits a quasi-complex orientation if and
only if, over each component~$B'$ of $B_\R \sminus Z$, the total
number of ovals of~$C_\R$ and
points of intersection of~$C_\R$ with the section~$0$
(see~\ref{real-trigonal})
equals $\chi(B') \bmod 2$.
\endRemark

\proposition\label{prop.orientation}
Any complex orientation of a non-hyperbolic trigonal curve
of type~$\I$ is quasi-complex.
\endproposition

\proof
For the extension of the orientation to~$\bar L$ one can mimic
the proof found in~\cite{F}. New is the case of a component
of~$B_\R$ that lies entirely in~$\Cal B$. The orientability
of~$\bar L$ over such a component
follows from Lemma~\ref{no.bold}.

For a non-hyperbolic trigonal curve~$C$ of type~$\I$, the
complement $C\sminus\pi_C^{-1}(B_\R)$ splits into four `quoters'
$C^\pm_\pm=C^\pm\cap C_\pm$.
Since both~$C_+$ and~$C_-$ are $c$-invariant,
whereas~$C^+$ and~$C^-$ are interchanged by~$c$,
any point of $C_\R$ over~$\Cal B$ is in the common part
of the boundaries $\partial C^\pm_-$. Thus, assuming
that $\pi_C(C^+_-) = B^+$, one concludes that
the map $\pi_C: C_\R \to B_\R$ is of degree~$+1$.
\endproof

\proposition\label{prop.hyperbolic}
Any hyperbolic trigonal curve~$C$ is of type~$\Ib$. Such a curve
is of type~$\I$ if and only if its base~$B$ is of
type~$\I$.
In
this case, with respect to some complex orientations of~$C_\R$
and~$B_\R$, one has $(\pi_C)_*[C_\R]=3[B_\R]$.
\endproposition

\proof
The first statement is a tautology, and the second one follows
immediately
from Lemma~\ref{base.type.I}. For the third statement, it suffices
to observe that, in the hyperbolic case, the halves~$C^\pm$ are
the pull-backs $\pi_C^{-1}(B^\pm)$.
\endproof

\section{Skeletons}\label{S.skeletons}

Unramified dessins of type~$\I$ can be reduced
to somewhat simpler objects, the so called skeletons,
which are obtained by disregarding all but dotted edges.
The principal result of this section
is Theorem~\ref{cor.Sk} describing
maximally inflected trigonal curves of type~$\I$
in terms of skeletons.

Throughout this section, we assume
that the underlying surface~$\base$ is orientable,
in other words, the base~$B$ of the ruling is of type~$\I$
and, hence, trigonal curves of type~$\Ib$ are those
of type~$\I$.

\subsection{Abstract skeletons}\label{a.skeletons}
Consider an embedded (finite) graph $\Sk \subset \bbase$ in a
compact surface $\bbase$. We do not exclude the possibility that
some of the vertices of~$\Sk$ belong to the boundary of~$\bbase$;
such vertices are called
\emph{real}. The set of edges at each real (respectively, inner)
vertex~$v$ of~$\Sk$ inherits from~$\bbase$ a pair of opposite linear
(respectively, cyclic) orders.
The \emph{immediate neighbors} of an edge~$e$ at~$v$ are the immediate
predecessor and successor of~$e$ with respect to (any) one of these orders.
A
\emph{first neighbor path} in~$\Sk$ is a sequence of
oriented
edges of~$\Sk$
such that each edge is followed by one of
its immediate neighbors.

Below, we consider graphs with connected components
of two kinds: directed and undirected. We call
such graphs \emph{partially directed}.
The directed and undirected parts of a partially directed
graph~$\Sk$ are denoted by~$\Skdir$ and~$\Skud$, respectively.
Accordingly,
we speak about directed and undirected vertices of these graphs.

\definition\label{def.a.skeleton}
Let~$\bbase$ be a compact
orientable
surface
with nonempty boundary.
An \emph{abstract skeleton} is a partially directed
embedded graph $\Sk \subset \bbase$,
disjoint
from the boundary $\partial\bbase$ except for some
vertices,
and satisfying the following conditions:
\roster
\item\local{Sk.2}
at each vertex of~$\Skdir$, the directions of adjacent edges
alternate;
\item\local{Sk.5} each real directed vertex has odd valency,
thus being a \emph{source} or a \emph{sink};
\item\local{Sk.6} each source is monovalent;
\item\local{Sk.1}
the graph~$\Sk$ has no first neighbor cycles and no inner isolated
vertices;
\item\local{Sk.3}
each boundary component~$l$ of~$\bbase$ has a vertex of~$\Sk$
and is subject
to the \emph{parity rule}\rom: directed and undirected vertices
alternate along~$l$;
\item\local{Sk.4}
if a component~$R$ of the cut $\bbase_{\Sk}$
has a single source in the boundary $\partial R$,
then~$R$ is a disk.
\endroster
\par\removelastskip
\enddefinition

\subsection{Admissible orientations}

\definition\label{def.adm.or}
Let~$\Sk \subset \bbase$ be an abstract skeleton.
An orientation of~$\Skud$
is called \emph{admissible} if, at each vertex, no two incoming
edges are immediate neighbors.
An \emph{elementary inversion} of an admissible orientation
is the
reversal
of the direction of one of the edges
so that the result is again an admissible orientation.
\enddefinition

\proposition\label{prop.adm.or1}
Any abstract skeleton~$\Sk$ has an admissible orientation.
Any two admissible orientations of~$\Sk$
can be connected by a sequence of elementary inversions.
\endproposition

This statement is proved
after Proposition~\ref{prop.adm.or2} below.
Due to the existence of admissible orientations,
all undirected edges of an abstract skeleton $\Sk \subset \bbase$
can be divided into two groups:
\emph{triggers}
and \emph{diodes}, the latter being those that
have the same direction, called \emph{natural},
in all admissible orientations of~$\Sk$.
On the contrary, each trigger has two \emph{states}
(orientations); each state can be extended to an admissible
orientation of the skeleton.

\proposition\label{prop.adm.or2}
Let $\Sk \subset \bbase$ be an abstract skeleton
and $e_1$, $e_2$ two distinct triggers.
Then, out of the four states of the pair $e_1$, $e_2$,
at least three extend to an admissible orientation.
\endproposition

\proof[Proof of Propositions~\ref{prop.adm.or1}
and~\ref{prop.adm.or2}]
To construct an admissible orientation of~$\Sk$, choose
an undirected edge~$e_1$, orient it arbitrarily,
and call the result $\vector{e}_1$ the \emph{first anchor}.
For each first neighbor path starting at $\vector{e}_1$,
orient each edge~$e'$ of the path in the direction
from $\vector{e}_1$, \ie, assign to~$e'$ the orientation induced
from the path. If the partial orientation thus obtained is
consistent, keep it; otherwise, disregard this orientation
and repeat the procedure starting from the first anchor $-\vector{e}_1$,
\ie, the same edge~$e_1$ with the orientation opposite to
the originally chosen one.

We assert that at least one of the anchors~$\vector{e}_1$
and~$-\vector{e}_1$ results in a consistent partial orientation.
Indeed, otherwise there are two oriented edges~$e'$, $e''$, a pair
of first neighbor paths $\gamma'_\pm$ starting
at~$\vector e_1$ and ending at $\pm e'$,
and a pair of first neighbor paths $\gamma''_\pm$ starting at
$-\vector e_1$ and ending
at $\pm e''$. Then, the sequence $\gamma'_+$, $\gamma^{\prime-1}_-$,
$\gamma''_+$, $\gamma^{\prime\prime-1}_-$ gives rise to a first
neighbor cycle.

If there is an edge~$e_2$ that has not yet
been assigned
an orientation, orient it and repeat the procedure starting
from the \emph{second anchor} $\vector{e}_2$ or $-\vector{e}_2$,
whichever works. Clearly, the new orientation is consistent with
the one obtained in the previous step. Continue this procedure
(selecting anchors and extending their orientation)
until the whole undirected part of~$\Sk$ is exhausted.

Notice that any admissible orientation of~$\Sk$ can be obtained
by this procedure with an appropriate choice of the anchors.

For the uniqueness up to elementary inversions, fix an admissible
orientation~$\frak o$ obtained from a sequence of anchors
$\vector{e}_1,\vector{e}_2,\ldots$ and consider another admissible
orientation~$\frak o'$. Among the first neighbor paths starting
at~$\vector{e}_1$ and ending at an edge~$e'$ whose $\frak
o'$-orientation differs from~$\frak o$, choose a maximal one. Then,
reversing the orientation of the last edge of such a maximal path is
an elementary inversion. Repeat this procedure
to switch to~$\frak o$ the orientations of all edges
reachable from $\vector{e}_1$; then continue with $\vector{e}_2$,
\etc.

To prove Proposition~\ref{prop.adm.or2},
note that an undirected edge~$e$ of $\Sk$
is a trigger if and only if both $\vector{e}$ and $-\vector{e}$
can be chosen for an anchor.
Choose a state~$\vector{e}_1$ of~$e_1$. If~$e_2$ can be reached
(by a first neighbor path) stating from~$\vector{e}_1$,
then $e_2$ cannot be reached from $-\vector{e}_1$,
as otherwise the two orientations induced
on~$e_2$ would coincide (no first neighbor cycles) and hence
$e_2$ would be a diode. Thus, one can assume
that~$e_2$ is not reachable from~$\vector{e}_1$,
and three admissible orientations
can be constructed
starting from the
anchors $\vector{e}_1,\pm\vector{e}_2,\ldots$
and $-\vector{e}_1,\ldots$.
\endproof

\midinsert
\eps{Fig-diod}
\figure\label{fig.diod}
A diod
\endfigure
\endinsert

\Remark\label{remark1}
If an abstract skeleton $\Sk \subset \bbase$ does not have cycles,
then any undirected edge of~$\Sk$ is a trigger.
Note however that in general diods do exist, see, \eg,
the fragment shown in Figure~\ref{fig.diod}. (This fragment can
easily be completed to a skeleton without inner vertices.)
Alternatively, one can consider a skeleton
with a monovalent inner vertex~$v$: the
only
edge adjacent to~$v$ is its own immediate neighbor;
hence it cannot be oriented towards~$v$.
\endRemark

\subsection{Equivalence of abstract skeletons}\label{equivalence}
Two abstract skeletons
are called \emph{equivalent} if, after a homeomorphism
of underlying surfaces,
they can be connected by a finite sequence of isotopies
and the following \emph{elementary moves},
\cf.~\ref{ss.equivalence}:
\Dashes
\dash \emph{elementary modification}, see Figure~\ref{fig.Sk};
\dash \emph{creating} (\emph{destroying}) \emph{a bridge}
in $\Skud$,
see Figure~\ref{fig.Sk};
the vertex shown in the figure can be inner or real,
and the dotted lines represent other edges of $\Skud$
that may
be present.
\endDashes
(A move is valid only if the result is again
an abstract skeleton.)
It is understood that
an
elementary move
does not mix~$\Skdir$ and~$\Skud$
and, when
acting on~$\Skdir$, a move must respect the prescribed
orientations (as shown in the figures), thus defining an
orientation on the resulting directed part.
On~$\Skud$, a move
is required to
respect some admissible orientation
of the original skeleton
and take
it to an admissible orientation of the result.
Note that, in a contrast to the definition of equivalence of
dessins, see~\ref{ss.equivalence},
respecting a certain orientation, either prescribed or
admissible, is an extra requirement here, \cf.
Remark~\ref{rem.orientation}.

\midinsert
\eps{Fig-Sk}
\figure\label{fig.Sk}
Elementary moves of skeletons
\endfigure
\endinsert

An equivalence of two abstract skeletons in the same surface
and with the same set of vertices is called \emph{restricted}
if the homeomorphism is identical and
the isotopies above can be chosen identical on the vertices.

\Remark\label{remark2}
For~$\Skud$, the orientation condition
in the definition above is restrictive only if
all edges involved are diodes; otherwise,
a required admissible orientation does exist
due to Proposition~\ref{prop.adm.or2}. For creating a bridge,
it would even suffice to assume that at least one of the immediate
neighbors of the bridge created is a trigger.
In particular, if~$\bbase$ is a disk, then any elementary
move is allowed, see Remark~\ref{remark1}.
\endRemark

\subsection{Dotted skeletons}\label{s.skeletons}
From now on, for simplicity, we confine ourselves to dessins
without anti-hyperbolic components. All such components
would be monochrome, see Proposition~\ref{no.monochrome},
and thus could easily been incorporated
using the concept of (partial) reduction, see~\cite{DIK}.

\paragraph\label{dotted.skeleton}
Intuitively, the dotted skeleton is obtained from a
dessin $\Gamma$ by disregarding all but \dotted- edges
and patching the latter through all \white-vertices.
According to Theorem~\ref{th.summary}, each inner
dotted edge of type~$\2$ retains a well defined orientation,
whereas an edge of type~$\3$ may be broken by \white-vertex,
and for this reason, its orientation may not be defined.
As types do not mix, edges and
pillars
of types~$\2$
and~$3$ would form separate components of the skeleton.

\definition\label{def.skeleton}
Let
$\Gamma \subset \base$
be an unramified dessin of type~$\I$
without anti-hyperbolic components.
The \emph{\rom(dotted\rom) skeleton} of~$\Gamma$
is the partially directed graph
$\Sk=\Sk_\Gamma\subset\bbase$
obtained from~$\Gamma$ as follows:
\Dashes
\dash contract each
pillar
to a single
point and declare this point a vertex of~$\Sk$;
\dash
patch each inner \dotted- edge
through its
\white-vertex, if there is one,
and declare the result an edge~$\Sk$;
\dash let~$\Skdir$ and~$\Skud$ be the images
of the edges and pillars of type~$\2$ and~$\3$, respectively,
each edge of type~$\2$ inheriting its
orientation from~$\Gamma$.
\endDashes
Here, $\bbase$ is the
surface
obtained from~$\base$ by contracting
each pillar to a single point.
\enddefinition

\proposition\label{prop.Sk=Sk}
The skeleton~$\Sk$ of a dessin~$\Gamma$
as in Definition~\ref{def.skeleton}
is an abstract skeleton in the sense of
Definition~\ref{def.a.skeleton}.
\endproposition

\proof
Properties~\iref{def.a.skeleton}{Sk.2}--\ditto{Sk.6}
follow immediately from
Theorem~\ref{th.summary}.
Each component
of $\partial\bbase$ has a vertex of~$\Sk$ due to our assumption
that $\Gamma$ has no anti-hyperbolic component,
and the parity rule
in~\iref{def.a.skeleton}{Sk.3} is a consequence
of Theorem~\iref{th.summary}{3}.
Property~\iref{def.a.skeleton}{Sk.4} is merely a restatement
of requirement~\itemref{tg}{tg-triangle} in the definition of
trichotomic graph, see Subsection~\ref{tg}.

For~\iref{def.a.skeleton}{Sk.1}, apply a sequence of
\white-outs along type~$\3$ dotted edges to convert~$\Gamma$ to
a dessin~$\Gamma'$ with the same skeleton~$\Sk$
and all \white-vertices real.
The orientation of dotted edges of~$\Gamma'$ induces
the prescribed orientation of~$\Skdir$ and
an admissible orientation of~$\Skud$, the \white-vertices
of~$\Gamma'$ residing in the real dotted edges connecting
outgoing inner dotted edges and/or \cross-vertices.
Thus, a first
neighbor cycle of~$\Sk$ would give rise to an \emph{oriented}
dotted cycle of~$\Gamma'$, which contradicts to the definition of
dessin, see~\ref{sketches}.
Finally, the inner vertices of~$\Sk$ are the images of
hyperbolic components of~$\Gamma$, which are necessarily adjacent
to inner dotted edges.
\endproof

\proposition\label{Sk.extension}
Any abstract skeleton $\Sk\subset\bbase$ is the
skeleton
of a certain dessin~$\Gamma$
as in Definition~\ref{def.skeleton}\rom;
any two such dessins can be connected by a sequence
of isotopies and elementary moves,
see~\ref{ss.equivalence},
preserving the skeleton.
\endproposition

\proposition\label{prop.Sk}
Let $\Gamma_1, \Gamma_2 \subset \base$
be two
dessins as in Definition~\ref{def.skeleton}\rom;
assume that~$\Gamma_1$ and~$\Gamma_2$
have the same pillars.
Then, $\Gamma_1$ and $\Gamma_2$
are related by a restricted equivalence
if and only if so are the corresponding
skeletons~$\Sk_1$ and~$\Sk_2$.
\endproposition

Propositions~\ref{Sk.extension} and~\ref{prop.Sk}
are proved in Subsections~\ref{proof.prop} and~\ref{proof.prop1}.
Here, we state the following immediate consequence.

\theorem\label{cor.Sk}
There is a canonical
bijection
between the set of
equivariant fiberwise deformation classes of maximally inflected
type~$\I$ real trigonal curves without
anti-hyperbolic components and the set of equivalence
classes
of abstract skeletons.
\qed
\endtheorem

\subsection{Proof of Proposition~\ref{Sk.extension}}\label{proof.prop}
The underlying surface~$\base$ containing~$\Gamma$
is the orientable blow-up of~$\bbase$ at the vertices of~$\Sk$:
each inner (boundary) vertex~$v$ is replaced with
the circle (respectively, segment) of directions at~$v$.
The circles and segments inserted are the
pillars
of~$\Gamma$.
Each source gives rise to a jump and is decorated accordingly;
all other pillars consist of dotted edges (the \white-vertices
are to be inserted later, see~\ref{ss.white})
with \cross-vertices at the ends.
The proper transforms of the edges of~$\Sk$
are the inner dotted edges of~$\Gamma$.

\paragraph\label{ss.white}
The blow-up produces
a certain
\dotted- subgraph $\Sk'\subset\base$.
Choose an admissible orientation of~$\Sk$, see
Proposition~\ref{prop.adm.or1}, regard it as an orientation of the
inner edges of~$\Sk'$, and insert a \white-vertex at the center of
each real dotted segment connecting a pair of outgoing inner edges
and/or \cross-vertices.

\paragraph\label{ss.solid}
Let~$\bar U \subset \bbase$ be a closed regular neighborhood
of~$\Skdir$ disjoint from~$\Skud$,
and~$U \subset \base$ be the preimage of~$\bar U$.
Shrink $U$ along $\partial\base$
so that the boundary $\partial U$ contains the \black-vertices
and
take for the inner solid edges of~$\Gamma$
the connected components of the inner part
of $\partial U$,
defining real solid edges and monochrome vertices accordingly.
Note that,
in view of Condition~\iref{def.a.skeleton}{Sk.1},
each connected component of the inner part
$\partial U\sminus\partial\base$ is an interval rather than a
circle, as a circle in $\partial U$ disjoint from $\partial\base$
would contract to a first neighbor cycle in~$\Skdir$.

\paragraph\label{ss.bold}
At this point, the closure of the complement $\base\sminus U$ should
be the union of the type~$\3$ regions of the dessin in question, and
the cut $(D \sminus U)_{\Sk'}$ should contain inner bold edges only.
Let~$R$ be a region of this cut. The \black-- and \white-vertices
of~$\Gamma$ define germs of bold edges at the boundary $\partial R$;
due to the parity rule~\iref{def.a.skeleton}{Sk.3}, incoming and
outgoing bold edges alternate along each component of $\partial R$.
Consider a disk $B^2$ with a distinguished oriented diameter~$d$ and
let $\varphi\:\partial R\to\partial B^2$ be an orientation
preserving covering taking the incoming/outgoing bold edges to the
corresponding points of $d\cap\partial B^2$. In view
of~\iref{def.a.skeleton}{Sk.4}, the map $\varphi$ extends to a
ramified covering $\bar\varphi\:R\to B^2$, which can be assumed
regular over~$d$, and it suffices to take for the inner bold edges
of~$\Gamma$ the components of the pull-back $\bar\varphi^{-1}(d)$.
This completes the construction of a dessin extending~$\Sk'$.

\paragraph\label{ss.uniqueness}
For the uniqueness, first observe that
a decoration of~$\Sk'$ with \white-vertices
is unique up to
isotopy and \white-ins/\white-outs along \dotted- edges.
Indeed, assuming
all \white-vertices real,
each such decoration is obtained from a
certain admissible orientation, see~\ref{ss.white},
which is unique up to a sequence
of elementary inversions, see
Proposition~\ref{prop.adm.or1}, and an elementary
inversion
results in a \white-in
followed by a \white-out at the other end of the edge reversed.
Thus, the distribution of the \white-vertices can be assumed
`standard'.

\paragraph\label{ss.1.solid}
The
union of the solid edges of any dessin~$\Gamma$
extending $\Sk'$ is the inner part of the boundary of
the shrunk preimage~$U \subset \base$ of a certain
closed
neighborhood~$\bar U$ of~$\Skdir$,
\cf.~\ref{ss.solid}.
(This
neighborhood~$U$ is the union of the type~$\2$ regions
of~$\Gamma$, see Theorem~\ref{th.summary}.) We assert that, for
each~$\Gamma$, there is a decreasing family
of neighborhoods~$U_t$, $t\in[0,1]$,
$U_0=U$, $U_{t'}\subset U_{t''}$ for $t'>t''$, composed of
isotopies and finitely many elementary modifications of the
boundary and such that $U_1$ is a regular neighborhood,
see~\ref{ss.solid},
and, for each regular value $t\in[0,1]$, replacing the solid edges
of~$\Gamma$ with the inner part of $\partial U_t$ results in a valid
dessin. Indeed, each component of the cut $U_{\Sk'}$ is a disk with
holes and handles, and it can be simplified by the following
operations: \Dashes \dash first, `replant' each handle by two
monochrome modifications, see Figure~\ref{fig.handle}; \dash next,
eliminate each hole by a monochrome modification, see
Figure~\ref{fig.hole}; \dash finally, by a number of monochrome
modifications, cut the resulting surface into triangles,
see~\iref{tg}{tg-triangle}.
\endDashes
It is a routine to check that each monochrome modification used
can be chosen to involve a pair of distinct solid edges
(due to~\iref{tg}{tg-triangle}, unless the component
in question already is a triangle,
it has at least two solid edges in the boundary)
and that, under this assumption, each intermediate
trichotomic graph is a valid dessin.

\midinsert
\eps{Fig-handle}
\figure\label{fig.handle}
`Replanting' a handle
\endfigure
\endinsert

\midinsert
\eps{Fig-hole}
\figure\label{fig.hole}
Eliminating a hole
\endfigure
\endinsert

\paragraph\label{ss.1.bold}
Since any two regular
neighborhoods of~$\Skdir$
are isotopic, it remains to
consider two dessins that differ by bold edges only.
Any collection of inner bold edges of
a valid dessin is obtained by the
construction of~\ref{ss.bold}, from an appropriate ramified
covering
$\bar\varphi\:R\to B^2$ (\cf. the passage from a dessin to
a $j$-invariant in~\cite{DIK}). Since the restrictions
$\varphi\:\partial R\to\partial B^2$ corresponding to the two
dessins are homotopic, the extension is unique up to homotopy in
the class of ramified covering, see~\cite{Natanzon}.
Hence, the two dessins are related by a sequence of bold
modifications.
\qed

It is items~\ref{ss.bold} and~\ref{ss.1.bold} in the proof why we
had to exclude from the consideration the case of the base curve
of type~$\II$, \ie, nonorientable $\bbase$.

\subsection{Proof
of Proposition~\ref{prop.Sk}}\label{proof.prop1}
The `only if' part is obvious: an elementary move of a dessin
either leaves its skeleton intact or results in its
elementary modification; in the latter case, a pair of edges of
the same type is involved, \ie, either both directed (and then the
orientation is respected) or both undirected (and then some admissible
orientation is respected, see~\ref{ss.white}).

For the `if' part, consider the skeleton~$\Sk$
at the moment of
a modification.
It can be regarded as the skeleton of a
dessin with inner monochrome vertices allowed
(see admissible trichotomic graphs in~\cite{DIK}),
and,
repeating the proof of
Proposition~\ref{Sk.extension}, one can see that $\Sk$ does
indeed
extend to a certain dessin. The extension remains a valid
dessin~$\Gamma$ before the modification as well. Hence, due to the
uniqueness given by Proposition~\ref{Sk.extension}, one can assume
that the original dessin is~$\Gamma$, and then the modification of
the skeleton is merely an elementary modification of~$\Gamma$.

Destroying a bridge of a skeleton is the same as destroying
a bridge of the corresponding dessin, and the inverse operation
of creating a bridge extends to a dessin equivalent
to the original one due to the uniqueness given
by Proposition~\ref{Sk.extension}.
\qed

\section{The case of the rational base}\label{S.rational}

In this section we prove Theorems~\ref{th.existence}
and~\ref{th.main} and attempt a constructive description
of maximally inflected trigonal curves
of type~$\I$ in rational ruled surfaces.
Note that, in the settings of Theorems~\ref{th.existence}
and~\ref{th.main}, each skeleton~$\Sk$ is a forest
in the disk, and all vertices of~$\Sk$ are on the boundary.

\subsection{Proof
of Theorem~\ref{th.existence}}\label{proof.existence}
The `only if' part is given by Proposition~\ref{prop.orientation}.
For the `if' part, consider a regular neighborhood $V \subset B$
of~$B_\R$. Under the assumptions on the orientation,
the germ $C'=\pi_C^{-1}(V)$ is separated by~$C_\R$, \cf. \eg\
\cite{V, Proof of Theorem 1.3.A}. On the other hand, since the
covering $\pi_C\:C\to B$ is unramified over the two disks
$B\sminus B_\R$, the curve~$C$ is obtained from~$C'$ by
attaching six disks; hence it remains separated.
\qed

\subsection{Proof of Theorem~\ref{th.main}}\label{proof.main}
Consider a trigonal curve~$C'$ as in the statement,
and let $\Sk'$
be the associated skeleton in the disk $\bbase \simeq \Cp1\!/c$.
Destroying all bridges (see Remark~\ref{remark2}),
one can assume that the valency of
each vertex of~$\Sk'$ is at most one, \ie, each undirected
component of~$\Sk'$ is either an isolated vertex (zigzag) or a
single edge connecting two vertices (ovals). Ignore the zigzags:
their position is uniquely recovered by the parity
rule~\iref{def.a.skeleton}{Sk.3}. Then, $\Sk'$ turns into a
collections of disjoint chords in the disk~$\bbase$,
directed and undirected,
connecting points in~$\partial\bbase$ of three types: sources,
sinks, and undirected vertices.

Let~$\Sk''$ be the skeleton associated to
the other curve~$C''$ with the
same modifications as above. Identify the vertices of~$\Sk''$ with
those of~$\Sk'$ according to the homeomorphism of the real parts.
Let~$l$ be a shortest chord of~$\Sk'$, \ie, such that
one of the two arcs
constituting
$\partial\bbase \sminus l$ is free of vertices of the skeletons.
Performing, if necessary, an elementary modification, one can assume
that $l$ is also an edge of~$\Sk''$.
Remove the part of the disk cut off by~$l$ and proceed by
induction, ending up with a pair of empty skeletons, which are
obviously equivalent. Thus, $\Sk''$ is equivalent to~$\Sk'$, and
Propositions~\ref{Sk.extension} and~\ref{prop.Sk}
imply the theorem.
\qed

\Remark\label{remark3}
Theorem~\ref{th.main} does not extend to not maximally inflected
curves, even of type~$\I$; an example of two non-equivalent
$M$-curves with isotopic real parts is found in~\cite{DIK}.
Nor does the theorem extend to the case of a base of
positive genus: the two skeletons shown
in Figure~\ref{fig.baraban}
(where $\bbase$ is a cylinder)
are obviously not related by restricted equivalence, as the
orientations shown prohibit
any elementary modification.
One can easily construct pairs of skeletons with isomorphic real
parts and not related by any equivalence, restricted or not.
\endRemark

\midinsert
\eps{Fig-baraban}
\figure\label{fig.baraban}
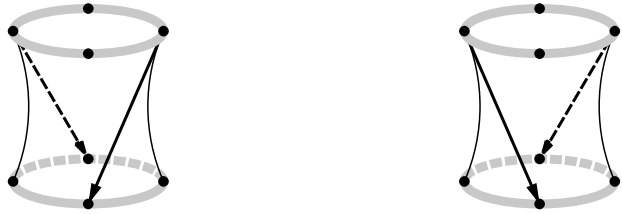
Nonequivalent skeletons with the same real part
\endfigure
\endinsert

\subsection{Blocks}\label{s.blocks}
In this section, we make an attempt of a constructive description
of the real parts of maximally inflected type~$\I$ trigonal curves
over the rational base.

\definition\label{def.block}
A type~$\I$ dessin~$\Gamma$ in the disk
is called a \emph{block} if~$\Gamma$
is unramified and
has no inner dotted edges of type~$\3$.
\enddefinition

It follows from Theorem~\ref{th.summary} that
all vertices of any block~$\Gamma$
are real,
all its ovals are of type~$\2$, and all
its zigzags are `short', \ie, each zigzag
contains a single \white-vertex.
In particular,
the real part of~$\Gamma$ consists of
$n=\frac13\deg\Gamma$ ovals and $n$ jumps, which are intermitted
with $2n$ zigzags; the position of the zigzags is uniquely
determined by the parity rule~\iref{th.summary}{3}.
Blocks are easily enumerated by the following statement.

\proposition\label{block.existence}
Let~$n\ge1$ be an integer, and let $O,J\subset S^1=\partial\base$ be
two disjoint sets of size~$n$ each. Then, there is a unique, up to
restricted equivalence,
block $\Gamma\subset\base$ of degree~$3n$
with an oval about each
point of~$O$, a jump at each point of~$J$, and a zigzag between
any two points of $O\cup J$ \rom(and no other pillars\rom).
\endproposition

\proof
Fix a bijection between~$J$ and~$O$ and connect each point of~$J$
to the corresponding point of~$O$ by a directed chord. Whenever
two chords intersect, resolve the crossing respecting the
orientation. Add to the resulting directed graph an isolated
vertex between any two points of $O\cup J$. The result is an
abstract dotted skeleton; due to Proposition~\ref{Sk.extension},
it extends to a block. The uniqueness is given by
Theorem~\ref{th.main}.
\endproof

Proposition~\ref{block.existence} and Theorem~\ref{th.summary}
provide a complete description of the real parts of maximally
inflected type~$\I$ trigonal curves
over~$\Cp1$.
Realizable are the real parts obtained as follows: start with a
disjoint union of a number of blocks, see
Proposition~\ref{block.existence}, and perform a sequence of
junctions converting the disjoint union of disks to a single
disk.

\Remark\label{remark4}
The description of maximally inflected
curves of type~$\I$ given above, in terms of
junctions, is similar to that of $M$-curves, see~\cite{DIK}.
However, unlike the case of $M$-curves,
in general a decomposition
of
an unramified
dessin of type~$\I$ into a junction of blocks is far
from unique.
\endRemark

\Remark\label{remark5}
Combining blocks, one can
also
obtain a great deal
of maximally inflected type~$\I$ trigonal curves
over irrational bases. However, in the case of the base
of positive genus, this construction is no longer universal:
there are unramified
dessins of type~$\I$ that cannot be cut into disks,
see, \eg, the skeletons
in Figure~\ref{fig.baraban}.
\endRemark

\appendix{}\label{ap.A}

In
the main part of the paper, we consider nonsingular trigonal curves
up to {\bf fiberwise} deformation equivalence,
\ie, we do not allow a pair of simplest (type $\tilde\bA_0^*$)
singular fibers to merge into a vertical flex.
This notion, natural in the framework of trigonal curves,
is not quite usual in general theory of nonsingular algebraic
curves in surfaces, where a less restrictive relation,
the so called rigid isotopy, is used. We reinterpret
this notion in terms of dessins and prove
that any non-hyperbolic nonsingular
real trigonal curve
of type~$\Ib$ is
rigidly isotopic to a maximally inflected one,
see Theorem~\ref{max.inflected}.

\subsection{Rigid isotopies and week equivalence}\label{s.rigid}
Keeping the conventional terminology,
we define \emph{rigid isotopy}
of nonsingular real trigonal curves
as the equivalence relation
generated by real isomorphisms and equivariant deformations
in the class of nonsingular (not necessarily almost-generic)
trigonal curves. Note that, in spite of the name `isotopy',
the underlying surface~$\Sigma$ and the base~$B$ are still not
assumed fixed: the complex structure is also subject
to deformation. Without this convention,
Proposition~\ref{equiv.zigzag} below would not hold.

Intuitively, the new notion differs from the deformation equivalence
by an extra pair of mutually inverse operations:
straightening/creating a zigzag, the former consisting in
bringing the two vertical tangents
bounding
a zigzag
together to a single vertical flex and pulling them apart
to the imaginary domain. On the level of dessins,
these operations are shown in Figure~\ref{fig.zigzag}.

\definition
Two dessins are called \emph{weakly equivalent} if they
are related by a sequence of isotopies,
elementary moves (see~\ref{ss.equivalence}),
and the operations
of \emph{straightening/\penalty0creating a zigzag}
consisting in replacing one of the fragments shown in
Figure~\ref{fig.zigzag} with the other one.
\enddefinition

\midinsert
\eps{Fig-zigzag}
\figure\label{fig.zigzag}
Straightening/creating a zigzag
\endfigure
\endinsert

The following statement is easily deduced from~\cite{DIK},
\cf.~Proposition~\ref{equiv.curves}.

\proposition\label{equiv.zigzag}
Two generic real trigonal curves are rigidly isotopic
if and only if their dessins are weakly equivalent.
\qed
\endproposition

\subsection{Creating zigzags}\label{creating-zigzags}
Let~$\Gamma$ be a dessin.
A \emph{lump} is a real component of~$\Gamma$
formed by two edges and two monochrome vertices
of the same kind.
Recall, see~\cite{DIK}, that $\Gamma$
is called \emph{bridge free} if
any bridge of~$\Gamma$ belongs to a lump.
(The two bridges forming a lump cannot be destroyed as this
operation would produce an oriented monochrome cycle.)
A \emph{long edge} in a bridge free dessin
is a sequence of inner edges interconnected
by lumps of the same kind.

A dessin is called \emph{peripheral} if it has no inner vertices
other than \cross-vertices.

In Statements~\ref{bold.not3}--\ref{zigzag+fat} below, $\Gamma$ is a
peripheral bridge free dessin of type~$\I$.
In particular, $\Gamma$ is not hyperbolic.

\lemma\label{bold.not3}
Let~$v$ be a real \white-vertex of~$\Gamma$ with the inner
bold edge of type~$\n3$. Then $v$ has a monochrome neighbor \rom(in
the boundary of the
adjacent
region~$R$ of type~$\1$\rom), followed by another
real \white-vertex~$u$ distinct from~$v$.
\endlemma

\proof
The only alternative is that the neighbor of~$v$ in
$\partial R$ is a real \cross-vertex. However, such a vertex
cannot be of type~$\1$, see Lemma~\ref{real.vertex}. Since $\Gamma$
is bridge free, a monochrome vertex must be followed by
another \white-vertex, switching the type of the real
\dotted-edges to~$\2$. Comparing these types shows that
$u\ne v$.
\endproof

\corollary\label{mono.not3}
Let $v$ be a \bold- monochrome vertex of~$\Gamma$ with the inner
edge of type~$\n3$. Then, up to equivalence of dessins, one can
assume that the real neighbors of~$v$ are \white-vertices.
\endcorollary

\proof
Assume that the real neighbors of~$v$ are \black-vertices and
consider the long inner edge~$e$ starting at~$v$. If the other
real end of~$e$ is monochrome, its real neighbors are two
\white-vertices. If the other end is a \white-vertex, then, due to
Lemma~\ref{bold.not3}, it has a monochrome neighbor followed by
another \white-vertex. In both cases, a \white-in followed by a
\white-out produces a desired dessin.
\endproof

\lemma\label{zigzag+fat}
Let~$v$ be an inner \cross-vertex of~$\Gamma$. Then $\Gamma$
is equivalent to a dessin
in which $v$ is included in a
fragment as in Figure~\ref{fig.zigzag}, right, possibly with
the dotted and/or bold inner edges long.
The new dessin
has at most one bridge.
\endlemma

\proof
Due to Lemma~\ref{inner.vertex}, the vertex~$v$ is necessarily of
type~$\1$ and, in view of Lemma~\ref{no13}, it must be connected
by a \solid-edge
(also of type~$\1$) to a
certain real \black-vertex~$u$. Using Corollary~\ref{mono.not3},
one can assume that the real neighbor connected to~$u$ by a real
\bold-edge is a \white-vertex. Then, performing, if necessary, a
\dotted-modification followed by
a \bold-modification or creating a \bold-bridge, one obtains a
desired
fragment as in Figure~\ref{fig.zigzag}, right.
\endproof

\proposition\label{to-unramified}
Any non-hyperbolic
dessin of type~$\I$ is weakly equivalent to an unramified
one.
\endproposition

\proof
Within a given weak equivalence class,
consider a dessin~$\Gamma$ with the minimal possible number of
inner \cross-vertices. According to~\cite{DIK}, one can assume~$\Gamma$
peripheral and bridge free.

Assume that $\Gamma$ has an inner \cross-vertex and show that it
can be taken out by creating a zigzag. Due to
Lemma~\ref{zigzag+fat}, $\Gamma$ can be replaced with a dessin
containing a fragment as in Figure~\ref{fig.zigzag}, right,
possibly with a number of lumps and edges long, and,
in order to create a zigzag, it remains to prove
that the lumps can be removed.

Let $u$, $v$, $w$, and~$m$ be the \black--, \cross--, \white--,
and monochrome vertices, respectively.

If $\Gamma$ has \bold- lump (in the long edge $[u,m]$),
proceed as follows.
\Dashes
\dash
If $m$ is part of a bridge, destroy the bridge and recreate it
back between $u$ and the first lump.
Otherwise, the next real neighbor~$w'$ of~$m$ is necessarily
a \white-vertex, and we examine the next real neighbor~$w''$
of~$w'$.
\dash
If $w''$ is monochrome, then a \white-vertex follows, and one
can apply a \white-in
and
destroy the bridge obtained and recreate it back as
above.
\dash
If $w''$ is a \black-vertex (necessarily $w''\ne u$ as
they are of different types), apply a
\white-in and a \black-in, slide the resulting inner \white-- and
\black-vertices through all but the last lumps, and perform
a \black-out and a \white-out \emph{to the last lump};
the resulting dessin has a desired fragment free
of \bold- lumps. (Note that the new fragment
is attached to another real component of the dessin.)
\endDashes

If $\Gamma$ has \dotted- lumps (in the long edge $[v,w]$),
proceed as follows.
\Dashes
\dash
If~$m$ is not part of a bridge (\ie, is followed by a
\white-vertex), then apply a \white-in, slide the inner \white-vertex
through the lumps, and apply a \white-out.
\dash
Otherwise, remove the bridge and consider
the long \bold- edge $[u,n]$ in the
`original' bridge free dessin.
By a \white-in, place an inner \white-vertex~$w'$ to this edge. (If $n$
is a \white-vertex, Lemma~\ref{bold.not3} is to be used.)
After a \dotted- monochrome modification, slide~$w'$ through
all \dotted- lumps
and, using the inverse modification and a \white-out
followed by creating a bridge, recreate~$\Gamma$ back,
now without \dotted- lumps.
\endDashes

Due to Lemma~\ref{no13},
the dessin~$\Gamma$
cannot have \solid- lumps. Hence, one can create a zigzag,
reducing the number of inner \cross-vertices of~$\Gamma$.
\endproof

\midinsert
\eps{Fig-star}
\figure\label{fig.star}
Essentially inner \cross-vertices
\endfigure
\endinsert

\theorem\label{max.inflected}
Any non-hyperbolic nonsingular
real trigonal curve
of type~$\Ib$ is
rigidly isotopic to a maximally inflected one.
\endtheorem

\proof
By a small equisingular perturbation one can make the curve
generic; then
the statement follows from Propositions~\ref{to-unramified}
and~\ref{equiv.zigzag}.
\endproof

\subsection{Further remarks}
In the disk, any non-hyperbolic
dessin (including those of type $\II$)
is also weakly equivalent to an unramified
one. This fact follows, \eg, from
Propositions~5.5.3 and~5.6.4 in~\cite{DIK},
see also~\cite{Zvonilov}.
Hence, any non-hyperbolic nonsingular real trigonal curve over the
rational base is rigidly isotopic to a maximally inflected
one.
In this statement, rigid isotopy can be understood in the
conventional sense, as an isotopy in the class of
nonsingular real algebraic
curves
in a fixed real ruled surface
$\Sigma\to\Cp1$.

The above statement does not extend directly to curves
over arbitrary bases.

\example\label{example}
Consider the dessin~$\Gamma$ shown in Figure~\ref{fig.star}. It is of
type~$\II$, and one can easily see that $\Gamma$ is not weakly
equivalent to an unramified dessin. Moreover, $\Gamma$ is
not equivalent to any dessin containing a fragment as in
Figure~\ref{fig.zigzag}, right, even with lumps. Indeed,
any such fragment would contain a \white-vertex, but all such
vertices are in odd hyperbolic components of~$\Gamma$,
one at each component,
and thus cannot be moved.
\endexample

\Remark\label{problem}
At present, we do not know whether two non-equivalent
maximally inflected trigonal curves of type~$\Ib$
can be rigidly isotopic. Note that this cannot happen
if the curves are $M$-curves, see~\cite{DIK}.
\endRemark

\refstyle{C}
\widestnumber\no{99}

\enddocument